\documentclass[11pt,a4paper]{amsart}

\usepackage{fullpage}
\usepackage{graphicx,amssymb}
\usepackage[top=35mm,bottom=35mm,left=30mm,right=30mm]{geometry}

\usepackage[frame]{xy}
\xyoption{arc}
\xyoption{graph}
\xyoption{matrix}

\theoremstyle{plain}
\newtheorem{theorem}{Theorem}[section]
\newtheorem{lemma}[theorem]{Lemma}
\newtheorem{corollary}[theorem]{Corollary}
\newtheorem{proposition}[theorem]{Proposition}

\theoremstyle{definition}
\newtheorem{definition}[theorem]{Definition}
\newtheorem{example}[theorem]{Example}
\newtheorem{remark}[theorem]{Remark}
\newtheorem{question}[theorem]{Question}

\let\le\leqslant
\let\ge\geqslant

\def\supp{\operatorname{supp}}

\def\dist{\operatorname{dist}}

\def\apex{\operatorname{apex}}
\let\DC\rightleftharpoons 

\def\<{\left\langle\, }
\def\>{\,\right\rangle}
\def\({\left(}
\def\){\right)}
\def\[{\left[\,}
\def\]{\,\right]}

\long\def\temp#1{\par\medskip\noindent\fbox{\begin{minipage}{.98\linewidth}\small #1\end{minipage}}}
\long\def\revised#1{\par\noindent\mbox{}\hskip-1em
\fbox{\begin{minipage}{1.1\linewidth}\footnotesize #1\end{minipage}}\par}

\begin{document}

\title{Uniqueness of quasi-roots \\
in right-angled Artin Groups}

\author{Eon-Kyung Lee}
\address{Department of Mathematics and Statistics, Sejong University, Seoul, Korea}
\email{eonkyung@sejong.ac.kr}

\author{Sang-Jin Lee}
\address{Department of Mathematics, Konkuk University, Seoul, Korea}
\email{sangjin@konkuk.ac.kr}


\begin{abstract}
We introduce the notion of quasi-roots and study their uniqueness
in right-angled Artin groups.
\end{abstract}

\maketitle


\section{Introduction}

In a group, if the equation $g^n=h$ holds,
then $g$ is called an \emph{$n$th root} of $h$.

The uniqueness of roots in groups has been studied by many authors,
e.g.~\cite{Bau60}.
Mal'cev proved that roots are unique (i.e.\ $g_1^n=g_2^n$ implies $g_1=g_2$)
in a torsion-free locally nilpotent group~\cite{Mal49}.
It is well-known that centralizers of
nontrivial elements of a torsion-free hyperbolic group
are necessarily infinite cyclic~\cite{BH99,GH90},
hence roots are unique.
In the Artin groups of finite type $\mathbf A_n$, ${\mathbf B}_n={\mathbf C}_n$
and affine type $\tilde{\mathbf A}_{n-1}$, $\tilde{\mathbf C}_{n-1}$,
roots are unique up to conjugacy~\cite{GM03,LL10,Mak71}.
For the uniqueness of roots in general mapping class groups, see~\cite{BP09}.
Pure braid groups and right-angled Artin groups are biorderable, hence
roots are unique~\cite{DT92,KR03,RZ98}.

In this article, we introduce the notion of quasi-roots
and study the uniqueness of quasi-roots in right-angled Artin groups.

\subsection{Right-angled Artin groups}
Before stating our results,
we recall basic notions in right-angled Artin groups.
See~\cite{Cha07, KK13, Ser89} for further details.

Throughout this article all graphs are finite and simple.
For a graph $\Gamma$, let $V(\Gamma)$ and $E(\Gamma)$ denote the vertex set
and the edge set of $\Gamma$, respectively.
For a graph $\Gamma$,
the \emph{right-angled Artin group} $A(\Gamma)$ on $\Gamma$ is defined by the presentation
$$ A(\Gamma)=\langle\, v\in V(\Gamma)\mid [v_i,v_j]=1\
\mbox{if $\{v_i,v_j\}\in E(\Gamma)$}\,\rangle.$$
In this article, we use the opposite convention
$$ G(\Gamma)=\langle\, v\in V(\Gamma)\mid [v_i,v_j]=1\
\mbox{if $\{v_i,v_j\}\not\in E(\Gamma)$}\,\rangle. $$
In other words, $G(\Gamma)=A(\Gamma^c)$, where $\Gamma^c$ denotes the complement graph of $\Gamma$.

Each element in $V(\Gamma)\cup V(\Gamma)^{-1}$ is called a \emph{letter}.
A \emph{word} means a finite sequence of letters.
For two words $w_1$ and $w_2$, the notation $w_1\equiv w_2$ means that
$w_1$ and $w_2$ coincide as sequences of letters.
For example, if $w_1\equiv v_1v_1v_2v_2^{-1}$ and $w_2\equiv v_1v_1$,
where $v_1,v_2\in V(\Gamma)$, then $w_1\not\equiv w_2$.
A word $w'$ is called a {\em subword} of a word $w$
if $w \equiv w_1 w' w_2$ for possibly empty words $w_1$ and  $w_2$.

An element $g$ in $G(\Gamma)$ can be expressed as a word $w$.
The word $w$ is called \emph{reduced} if $w$ is a shortest
word among all words representing $g$.
In this case, the length of $w$ is called the \emph{word length} of $g$,
denoted by $|g|$.

The \emph{support} of $g\in G(\Gamma)$, denoted by $\supp(g)$,
is defined as the set of
all generators $v\in V(\Gamma)$ such that $v$ or $v^{-1}$ appears in a reduced word
representing $g$.
It is known that $\supp(g)$ is well-defined.
For example, see \cite{KK13,KK14,LL16,LL18}.

Let $w$ be a (possibly non-reduced) word in  $G(\Gamma)$.
A subword $v^{\pm1}w_1v^{\mp1}$, $v\in V(\Gamma)$, of $w$ is called a \emph{cancellation} of $v$ in $w$
if each generator in ${\operatorname{supp}}(w_1)$ commutes with $v$.
If, furthermore, no letter in $w_1$ is equal to $v$ or $v^{-1}$,
it is called an \emph{innermost cancellation} of $v$ in $w$.
It is known that $w$ is reduced if and only if $w$ has no innermost cancellation
~\cite{KK15,LL16,LL18}.

Abusing notation, we shall sometimes regard a word as the group element represented by that word.

\begin{definition}[disjointly commute]
We say that $g_1, g_2\in G(\Gamma)$ \emph{disjointly commute},
denoted by $g_1\DC g_2$,
if $\supp(g_1)\cap\supp(g_2)=\emptyset$ and  each $v_1\in \supp(g_1)$ commutes with each $v_2\in\supp(g_2)$
in $G(\Gamma)$, i.e.\  $\{v_1,v_2\}\not\in E(\Gamma)$.
\end{definition}

\begin{definition}[strongly non-split]
An element $g\in G(\Gamma)$ is \emph{non-split}
if $\supp(g)$ spans a connected subgraph of $\Gamma$.
Otherwise, it is \emph{split}.
We say that $g$ is \emph{strongly non-split} if $g$ is non-split
and there is no generator $v\in V(\Gamma)$ that disjointly commutes with $g$.
\end{definition}

\begin{example}\label{eg:1}
Let $P_5=${\small $\begin{xy}/r.6mm/:
(0,0)="a" *+!U{v_1} *{\bullet},
"a"+(10,0)="a" *+!U{v_2} *{\bullet},
"a"+(10,0)="a" *+!U{v_3} *{\bullet},
"a"+(10,0)="a" *+!U{v_4} *{\bullet},
"a"+(10,0)="a" *+!U{v_5} *{\bullet},
(0,0)="a"; "a"+(40,0) **@{-},
\end{xy}$}
be the path graph on five vertices
$v_1,\ldots,v_5$.
The right-angled Artin group $G(P_5)$ has a presentation
$$
G(P_5)=\langle v_1,\ldots,v_5\mid [v_i,v_j]=1~\mbox{if $|i-j|\ge 2$}\rangle.
$$
Let $g_1=v_2v_3 v_5^{-2}$, $g_2=v_2v_3$ and $g_3=v_2v_3v_4$.
Then $\supp(g_1)=\{v_2,v_3,v_5\}$,
$\supp(g_2)=\{v_2,v_3\}$ and
$\supp(g_3)=\{v_2,v_3,v_4\}$.
The element $g_1$ is split because $\supp(g_1)$
spans the disconnected subgraph
{\small $\begin{xy}/r.6mm/:
(0,2)="org",
"org"="a" *+!U{v_2} *{\bullet},
"a"+(10,0)="a" *+!U{v_3} *{\bullet},
"a"+(10,0)="a" *+!U{v_5} *{\bullet},
"org"="a"; "a"+(10,0) **@{-},
\end{xy}$}
 of $\Gamma$.
The element $g_2$ is non-split but it is not strongly non-split
because $v_5$ disjontly commutes with $g_2$.
The element $g_3$ is strongly non-split
because $\supp(g)=\{v_2,v_3,v_4\}$ spans the connected subgraph
{\small $\begin{xy}/r.6mm/:
(0,2)="org",
"org"="a" *+!U{v_2} *{\bullet},
"a"+(10,0)="a" *+!U{v_3} *{\bullet},
"a"+(10,0)="a" *+!U{v_4} *{\bullet},
"org"="a"; "a"+(20,0) **@{-},
\end{xy}$}
of $\Gamma$
and because no generator disjointly commutes with $g$.
\end{example}

\begin{definition}[geodesic decomposition]
For elements $g,g_1,\ldots,g_k\in G(\Gamma)$, the decomposition
$g=g_1\cdots g_k$ is a \emph{geodesic decomposition} of $g$ if
$|g|=|g_1|+\cdots+|g_k|$.
\end{definition}

\begin{definition}[cyclic conjugate]
For $g\in G(\Gamma)$ with a geodesic decomposition $g=g_1g_2$,
the element $g_2g_1$ is a {\em cyclic conjugate} of $g$.
For a word $w\equiv w_1 w_2$, the word $w_2 w_1$
is  a {\em cyclic conjugate} of $w$.
\end{definition}

\begin{definition}[cyclically reduced]
An element $g\in  G(\Gamma)$ is {\em cyclically reduced}
if $g$ has the minimum word length in its conjugacy class.
\end{definition}

Servatius~\cite[Proposition on p.38]{Ser89} proved that for any $g\in G(\Gamma)$
there exists a geodesic decomposition $g=u^{-1}hu$
for $u,h\in G(\Gamma)$ such that $h$ is cyclically reduced.

\begin{definition}[primitive]
A nontrivial element $g\in G(\Gamma)$ is \emph{primitive}
if $g$ is not a nontrivial power of another element,
i.e.\ $g=u^n$ never holds for any $n\ge 2$ and $u\in G(\Gamma)$.
Similarly, a nonempty word $w$ is \emph{primitive}
if $w$ is not a nontrivial power of another word.
\end{definition}

\subsection{Quasi-roots}
We define quasi-roots as follows.

\begin{definition}
Let $0\le \lambda < \frac12$ and $N\ge 2$.
For $h\in G(\Gamma)$, an element
$g\in G(\Gamma)$ is a \emph{$(\lambda,N)$-quasi-root} of $h$
if there exists a geodesic decomposition
\begin{equation}\label{E:qr}
h=ag^n b=a\cdot \underbrace{g\cdot g\cdots g}_n\cdot b
\end{equation}
for some $n\ge N$ and $a,b\in G(\Gamma)$ with $|a|, |b|\le\lambda|h|$.
The above decomposition is called a {\em $(\lambda,N)$-quasi-root decomposition} of $h$.
\end{definition}

We are interested in the uniqueness of quasi-roots:
given an element $h\in G(\Gamma)$,
we want to find conditions under which we can guarantee
that a quasi-root of $h$ is uniquely determined up to conjugacy.

Here are immediate observations concerning the quasi-root decomposition \eqref{E:qr}.
\begin{enumerate}

\item
If $|a|=|b|=0$, then \eqref{E:qr} becomes $h=g^n$,
hence $g$ is an $n$th root of $h$.

\item
The condition $N\ge 2$ is natural
because  we cannot expect any kind of
uniqueness of quasi-roots
for $n=1$ (and $|a|+|b|>0$).
Meanwhile, since $n\ge N\ge 2$ and $g^n=g\cdots g$ is geodesic,
$g$ is cyclically reduced. (See Lemma~\ref{lem:c-red1}.)

\item
The condition $\lambda<\frac12$ is natural.
If $|a|\ge \frac12|h|$ or $|b|\ge \frac12|h|$ is allowed,
then $g$ is not uniquely determined up to conjugacy.
For example, let $h=v_1^nv_2^n$, where $v_1$ and $v_2$ are different generators.
Then $h=v_1^n\cdot 1^n\cdot v_2^n=1\cdot v_1^n\cdot v_2^n
= v_1^n\cdot v_2^n\cdot 1$,
hence $g$ is possibly one of  $1$, $v_1$ and $v_2$.
However, the elements $1$, $v_1$ and $v_2$ are not conjugate to each other.
\end{enumerate}

The following is the main result of this article.

\begin{theorem}\label{thm:AA}
Let $\Gamma$ be a finite connected graph, and
let $0\le \lambda < \frac12$ and $N\ge \frac{2|V(\Gamma)|+1}{1-2\lambda}$.
If there are two $(\lambda,N)$-quasi-root decompositions
$$h=a_1g_1^{n_1}b_1=a_2g_2^{n_2}b_2$$
such that $g_1$ and $g_2$ are strongly non-split and primitive,
then the quasi-roots $g_1$ and $g_2$ are conjugate such that
$a_1g_1a_1^{-1}=a_2g_2a_2^{-1}$ and $b_1^{-1}g_1b_1=b_2^{-1}g_2b_2$.
\end{theorem}

The constant $N$ depends only on $\lambda$ and $|V(\Gamma)|$.
Hence we obtain the following.

\begin{corollary}\label{thm:A}
Let $\Gamma$ be a finite connected graph.
For any $\lambda\in[0,\frac12)$, there exists $N\ge 3$ such that
strongly non-split and primitive $(\lambda,N)$-quasi-roots
are unique up to conjugacy.
\end{corollary}

Concerning the main theorem, we remark the following.
\begin{enumerate}
\item
The strongly non-splitness of $g_i$ is a necessary condition
for the uniqueness of quasi-roots.
For example, consider $G(P_5)$ in Example~\ref{eg:1}.
Let $h=(v_2^mv_3^m v_5)^n$.
Observe that $h=1\cdot (v_2^mv_3^m v_5)^n\cdot 1
=v_5^n\cdot (v_2^mv_3^m)^n\cdot 1
=1\cdot (v_2^mv_3^m)^n\cdot v_5^n$.
Therefore, if $m$ is large enough so that $(2m+1)\ge 1/\lambda$,
then $n\le \lambda(2m+1)n=\lambda|h|$.
Hence both  $v_2^mv_3^mv_5$ and $v_2^mv_3^m$ are $(\lambda,n)$-quasi-roots of $h$.
However, $v_2^m v_3^m v_5$ and $v_2^mv_3^m$ are not conjugate.
Notice that $v_2^m v_3^m v_5$ is split and
that $v_2^mv_3^m$ is non-split but not strongly non-split.

\item
The primitiveness of $g_i$ is not a restriction.
Without that condition, one may conclude that
there exists an element $u$
such that any $(\lambda,N)$-quasi-root is conjugate to a power of $u$.
\end{enumerate}

Note that we cannot expect a naive generalization of the main theorem
to mapping class groups.
For example, consider the $n$-strand Artin braid group $B_n$.
Let $\delta=\sigma_{n-1}\cdots\sigma_1$ and
$\epsilon=\delta\sigma_1$,
where $\sigma_1,\ldots,\sigma_{n-1}$ are the standard generators of $B_n$.
Then $\delta^n=\epsilon^{n-1}$.
For any $m\ge 1$, the element $h=\delta^{mn}=\epsilon^{m(n-1)}$ has two primitive (quasi-)roots
$\delta$ and $\epsilon$ which are not conjugate to each other.
In mapping class groups, pseudo-Anosov elements
may behave like strongly non-split elements in right-angled Artin groups.
Hence the following questions would be interesting.

\begin{question}
In mapping class groups,
is it true that primitive pseudo-Anosov quasi-roots are unique up to conjugacy?
\end{question}

\begin{question}
In torsion-free hyperbolic groups, is it true that primitive quasi-roots are unique up to conjugacy?
\end{question}

The results in this article seem insufficient to build a finite-time algorithm 
for finding a quasi-root or deciding the existence of a quasi-root.
Hence the following questions would be also interesting.

\begin{question} [Quasi-root decision problem] \mbox{}\\
Given a graph $\Gamma$ and $(\lambda,N, h)  \in  [0, 1/2) \times  [2, \infty)  \times  G(\Gamma)$, 
decide whether or not $h$ has a $(\lambda,N)$-quasi-root.
\end{question}

\begin{question}[Quasi-root search problem] \mbox{} \\
Given a graph $\Gamma$ and $(\lambda,N, h) \in [0, 1/2) \times  [2, \infty)  \times  G(\Gamma)$,
where $h$ is assumed to have a $(\lambda,N)$-quasi-root,
find a $(\lambda,N)$-quasi-root $g$ of $h$ together with a $(\lambda,N)$-quasi-root decomposition $h=ag^n b$.
\end{question}

\section{SD-conical elements}

Crisp, Godelle and Wiest introduced a normal form
of elements and the notion of pyramidal elements
in order to solve the conjugacy problem in right-angled Artin groups
in linear-time~\cite{CGW09}.
We use their normal form and a variation
of pyramidal elements, called SD-conical elements.
In this section we define SD-conical elements and explore
their properties.

The following equivalences are used implicitly in many works.
We include a proof for readers' convenience.

\begin{lemma}\label{lem:c-red1}
The following are equivalent for $g\in  G(\Gamma)$.
\begin{enumerate}
\item
$g$ is cyclically reduced.

\item There is no geodesic decomposition as $g=u^{-1}hu$
for $u, h\in G(\Gamma)$ with $u\ne 1$.

\item
$g^n = gg\cdots g$ is geodesic (i.e.\ $|g^n|=n|g|$) for all $n\ge 2$.

\item
$g^n = gg\cdots g$ is geodesic (i.e.\ $|g^n|=n|g|$) for some $n\ge 2$.

\item
Every cyclic conjugate of $g$ has the same word length as $g$,
i.e.\ for any geodesic decomposition $g_1 g_2$ of $g$,
the decomposition $g_2 g_1$ is also geodesic.
\end{enumerate}
\end{lemma}

\begin{proof}
``(i) $\Leftrightarrow$ (ii)'' is known
by Servatius \cite[Proposition on p.38]{Ser89}.
The implications ``(iii) $\Rightarrow$ (iv) $\Rightarrow$ (v) $\Rightarrow$ (ii)''
are obvious. Therefore it is enough to show ``(ii) $\Rightarrow$ (iii)''.

Assume that $g^n=gg\cdots g$ is not geodesic for some $n\ge 2$.
Let $w$ be a reduced word representing $g$.
Then $w^n=ww\cdots w$ has an innermost cancellation.
Since $w$ is reduced, this cancellation occurs between two consecutive $w$'s,
i.e.\ there are words $w_1,\ldots, w_4$ and a letter $x\in V(\Gamma)^{\pm 1}$
such that
$$
w\equiv w_1xw_2\equiv w_3x^{-1}w_4,\qquad
x\DC w_2,\qquad
x^{-1}\DC w_3.
$$
Since $x\DC w_2$, we have $x\not\in\supp(w_2)$,
hence the subword $w_3x^{-1}$ must be a prefix of $w_1$.
Therefore $w_1\equiv w_3x^{-1}w_5$ for some word $w_5$.
Now we have
$$
w\equiv w_3x^{-1}w_5xw_2,\qquad
x\DC w_2,\qquad
x^{-1}\DC w_3.
$$
Then $g$ has a geodesic decomposition as $g=x^{-1}hx$, where $h=w_3w_5w_2$.
This contradicts (ii).
\end{proof}

\begin{definition}[set of starting generators]
For $g\in G(\Gamma)$, the \emph{set of starting generators}, denoted by $S(g)$,
is the set of all generators $v\in V(\Gamma)$ such that either $v$ or $v^{-1}$ is
the first letter of a reduced word representing $g$,
i.e.
$$S(g)=\{v\in V(\Gamma)\mid g=v^\epsilon h
\mbox{~is geodesic for some $\epsilon=\pm 1$
and $h\in G(\Gamma)$}\}.$$
\end{definition}

\begin{definition}[conical element and apex]
For $g\in G(\Gamma)$, if $S(g)$ consists of a single generator, say $v_0$,
then we say that $g$ is \emph{conical} or \emph{$v_0$-conical}.
The generator $v_0$ is called the {\em apex} of $g$,
denoted by $\apex(g)$.
\end{definition}

Note that a conical element is necessarily non-split.

The following proposition is a quantitative version of Proposition 2.18 in~\cite{CGW09}.
Though it can be proved in almost the same way,
we include a sketchy proof for better readability.

\begin{proposition}\label{prop:pc}
Let $g \in G(\Gamma)$ be non-split and cyclically reduced,
and let $v_0\in\supp(g)$. Then $g$ is conjugate to a
$v_0$-conical element $p\in G(\Gamma)$ by $a,b\in G(\Gamma)$ as
$$g = apa^{-1} = b^{-1} p b,$$
where $g^k = ab$ holds and is geodesic for some $0\le k\le |V(\Gamma)|-1$.
In particular, if $n\ge |V(\Gamma)|$, then $g^n$ has a geodesic decomposition
as $g^n = a p^{n-k} b$.
\end{proposition}

\begin{proof}[Sketchy proof]
If $g$ is $v_0$-conical, there is nothing to prove because
we can take $k=0$ and $a=b=1$.
So we assume that $g$ is not $v_0$-conical.

We follow the argument in~\cite{CGW09}.
For any non-split element $h\in G(\Gamma)$ and for any $v\in\supp(h)$,
there exists a unique geodesic decomposition as $h=tp$
such that $p$ is $v$-conical and $v\not\in\supp(t)$~\cite[Lemma 2.16]{CGW09}.

Starting from $g_0=g$, we inductively do the following for $i=0,1,2,\ldots$:
(i) determine the geodesic decomposition $g_i=t_i\cdot p_i$,
where $p_i$ is $v_0$-conical and $v_0\not\in\supp(t_i)$;
(ii) take the conjugation $g_{i+1}=t_i^{-1}g_it_i=p_it_i$.

Then $g_j$ is $v_0$-conical for some $1\le j\le |V(\Gamma)|-1$.
(In fact, the proof of \cite[Proposition 2.18]{CGW09} says that
one can take $j\le \max\{\dist_{\Gamma_1}(v_0,v)\mid v\in\supp(g)\}$,
where $\Gamma_1$ is the subgraph of $\Gamma$ spanned by $\supp(g)$.
Notice that
$\max\{\dist_{\Gamma_1}(v_0,v)\mid v\in\supp(g)\}
\le |V(\Gamma)|-1$.)
Let $k$ be the smallest integer such that $g_k$ is $v_0$-conical.
Then $1\le k\le |V(\Gamma)|-1$.

In the decomposition $g_k=t_k\cdot p_k$, we have $t_k=1$ and $g_k=p_k$.
Keeping track of the above procedure, we get
\begin{align*}
g_0 & =  g = t_0 \cdot p_0, \\
g_1 & =   t_0^{-1}g_0 t_0 = p_0 t_0 = t_1 \cdot p_1, \\
g_2 & =   t_1^{-1}g_1 t_1 =p_1 t_1 = t_2 \cdot p_2, \\
&\cdots\\
g_k & =   t_{k-1}^{-1}g_{k-1}t_{k-1} = p_{k-1}t_{k-1} = t_kp_k=p_k.
\end{align*}

Let $p=g_k=p_k$, $a = t_0 \cdots t_{k-1}$ and $b = p_{k-1} \cdots p_0$.
Then the following hold
$$ |g|=|g_0|=\cdots = |g_{k-1}| = |g_k|=|p|,$$
$$|g_i^m| = m|g_i| \mbox{ for } 0\le i\le k,\ m\ge 1$$
because $g$ is cyclically reduced and each $g_{i+1}$ is a cyclic conjugate of $g_i$,
and hence also cyclically reduced for $0\le i\le k-1$.

Using the identities $p_i t_i = t_{i+1} p_{i+1}$ for $0\le i\le k-1$, we obtain
\begin{align*}
g^k &= (t_0 p_0)^k =  \underbrace{(t_0 p_0)  (t_0 p_0) \cdots (t_0 p_0)}_k =
 ( t_0 t_1 \cdots t_{k-1}) (p_{k-1} \cdots p_1 p_0) = ab, \\
ga & = t_0 p_0 t_0 t_1 \cdots t_{k-1} = t_0 \cdots t_{k-1} t_k p_k = a p, \\
bg &= p_{k-1} \cdots p_1 p_0 t_0 p_0 = t_k p_k p_{k-1} \cdots p_1 p_0 = pb,
\end{align*}
namely $g^k = ab$ and $g = apa^{-1} = b^{-1} p b$.

On the other hand,
$$|a|+|b| \leq |t_0|+\cdots + |t_{k-1}|+ |p_0|+\cdots + |p_{k-1}|
= |g_0|+\cdots + |g_{k-1}| = k|g| = |g^k|,$$
where the first equality holds because each $g_i=t_i\cdot p_i$ is a geodesic decomposition.
Since  $|g^k|\leq |a|+|b|$ is obvious, we get $|a|+|b|=|g^k|$, namely $g^k = ab$ is a geodesic decomposition.

The identity $g^n = a p^{n-k} b$ follows from $g^k = ab$ and $g =  b^{-1} p b$.
Observe
$$|a|+ | b| + | p^{n-k}|  = k|g| + (n-k)|p| = n|g|=|g^n|,$$
hence $g^n = a p^{n-k} b$ is a geodesic decomposition.
\end{proof}

\begin{definition}[$v_0$-conical conjugate]
Given a non-split and cyclically reduced element $g\in G(\Gamma)$
and a generator $v_0\in\supp(g)$,
the procedure in the proof of Proposition~\ref{prop:pc} uniquely determines
a $v_0$-conical element $p$, called
the \emph{$v_0$-conical conjugate} of $g$.
\end{definition}

The conical conjugate $p$
in the above definition
is also non-split and cyclically reduced
because $p$ is obtained from $g$
by iterated cyclic conjugations and
cyclic conjugations preserve non-splitness and cyclically reducedness.
In particular, $|g|=|p|$.

\begin{definition}[CGW-normal form]
Suppose $V(\Gamma)$ is endowed with a linear order.
\begin{enumerate}
\item
A reduced word $w\equiv v_{i_1}^{\epsilon_1}\cdots v_{i_k}^{\epsilon_k}$
on $V(\Gamma)^{\pm 1}$ representing an element $g\in G(\Gamma)$
is \emph{initially normal}
if either $w$ is trivial or
$v_{i_1}$ is the largest element of $S(g)$.

\item A reduced word $w$ is \emph{normal} if all its suffixes are initially normal.

\item
Any element $g\in G(\Gamma)$ has a unique normal representative
word~\cite[Proposition 2.6]{CGW09}.
We call this normal word the \emph{CGW-normal form} of $g$,
and denote it by $\sigma(g)$.
\end{enumerate}
\end{definition}

Note that all the subwords of a normal word are normal words.

One can understand $\sigma(g)$ as follows.
Extend the linear order on $V(\Gamma)$ to $V(\Gamma)^{\pm 1}$
so that each generator $v$ is the immediate predecessor of $v^{-1}$.
(Namely, if the order on $V(\Gamma)$ is $v_1<\cdots<v_r$, then
the order on $V(\Gamma)^{\pm 1}$ is
$v_1<v_1^{-1}<v_2<v_2^{-1}<\cdots <v_r<v_r^{-1}$.)
Then $\sigma(g)$ is the largest in the lexicographic order
among all the reduced words representing $g$.

\begin{definition}[SD-conical]
Suppose that $V(\Gamma)$ is endowed with a linear order $<$.

(i) A conical element $g$ is \emph{pyramidal} if $\apex(g)$ is the
smallest element of $\supp(g)$.

(ii) A conical element $g$ is \emph{SD-conical}
if $v_0=\apex(g)$ does not commute with
any generator smaller than $v_0$, i.e.\
if $v< v_0$, then $\{v,v_0\}\in E(\Gamma)$.
(From a graph theoretical point of view, the vertex $v_0$ dominates
all the vertices $v\in V(\Gamma)$ smaller than $v_0$.
SD stands for \emph{smaller vertex dominating}.)
\end{definition}

\begin{example}
Consider $G(P_5)$ in Example~\ref{eg:1}.
Define a linear order on $V(P_5)$ by $v_1<v_2<v_3<v_4<v_5$.
Let $g=v_2v_4^{-1}v_3^{-1}v_5$.
Then $S(g)=\{v_2,v_4\}$.
The following are reduced words representing $g$.
$$w_1\equiv v_2v_4^{-1}v_3^{-1}v_5,\quad
w_2\equiv v_4^{-1}v_2v_5v_3^{-1},\quad
w_3\equiv v_4^{-1}v_5v_2v_3^{-1},\ldots
$$
The word $w_1$ is not initially normal because $v_2$
is not the largest element of $S(g)$.
The word $w_2$ is initially normal, but it is not normal
because the suffix $v_2v_5v_3^{-1}$ is not initially normal.
The word $w_3$ is normal.
Hence $\sigma(g)\equiv w_3\equiv v_4^{-1}v_5v_2v_3^{-1}$.

The element $g=v_2v_4^{-1}v_3^{-1}v_5=v_4^{-1}v_5v_2v_3^{-1}$ is not conical
because $S(g)$ is not a singleton.
The conjugate $g_1=(v_4^{-1}v_5)^{-1}g(v_4^{-1}v_5)=v_2v_3^{-1}v_4^{-1}v_5$
is pyramidal because $S(g_1)=\{v_2\}$ and $v_2=\apex(g_1)$ is the
smallest element of $\supp(g_1)=\supp(g)=\{v_2,v_3,v_4,v_5\}$.
Moreover, $g_1$ is SD-conical because
$v_1$ is the only generator smaller than $v_2=\apex(g_1)$
but it does not commute with $v_2$.

In the same way, one can see that
$h_1=v_4v_5$ is pyramidal but not SD-conical,
and that $h_2=v_2v_1v_3v_4$ is SD-conical but not pyramidal.
\end{example}

The notions of normal form, pyramidal and SD-conical are
defined via a linear order on $V(\Gamma)$.
Whenever such terms are used,
it is assumed that $V(\Gamma)$ is endowed with
some linear order even though it is not specified.

\begin{lemma}\label{lem:nf}
Let $g_1g_2$ be geodesic for $g_1,g_2\in G(\Gamma)$.
\begin{enumerate}
\item
If $g_2$ is SD-conical,
then $\sigma(g_1g_2)\equiv \sigma(g_1)\sigma(g_2)$.

\item
If $g_1$ is SD-conical and strongly non-split,
then $g_1g_2$ is also SD-conical and strongly non-split.
\end{enumerate}
\end{lemma}

\begin{proof}
In the proof of Proposition 2.20 in~\cite{CGW09}, the following is observed.
\begin{itemize}\item[]
If $v_1^{\epsilon_1}, v_2^{\epsilon_2}$ are letters and $w$ is a reduced word
such that the words $v_1^{-\epsilon_1}w$ and $wv_2^{\epsilon_2}$
are reduced but the word $v_1^{-\epsilon_1}wv_2^{\epsilon_2}$ is not,
then $v_1^{\epsilon_1}=v_2^{\epsilon_2}$ and $v_1\DC w$.
\end{itemize}

Using this, we get the following.

\medskip
\noindent\textsl{Claim.}\ \
For $v\in V(\Gamma)$,
if $v\in S(g_1g_2)$ and $v\not\in S(g_1)$,
then $v\in S(g_2)$ and $v\DC g_1$.

\begin{proof}[Proof of Claim]
Let $w_1$ and $w_2$ be reduced words representing $g_1$ and $g_2$, respectively.
Since $v\in S(g_1g_2)$ and $v\not\in S(g_1)$, we can write
$w_2\equiv w_2'v_2^{\epsilon_2}w_2''$ such that
$v\not\in S(w_1w_2')$ and $v\in S(w_1w_2'v_2^{\epsilon_2})$, i.e,
$v^{-\epsilon}w_1w_2'$ is reduced but
$v^{-\epsilon}w_1w_2'v_2^{\epsilon_2}$ is not.
Since $g_1g_2$ is geodesic, $w_1w_2'v_2^{\epsilon_2}$ is reduced.
From the above observation,
$v\DC w_1w_2'$ and $v^{\epsilon}=v_2^{\epsilon_2}$.
Hence
$v\DC g_1$ and $v\in S(g_2)$.
\end{proof}

\bigskip
(i)\ \
Our proof is similar to the one in~\cite[Proposition 2.20]{CGW09},
where it is shown that if $p$ is pyramidal and cyclically reduced,
then $\sigma(p)\sigma(p)$ is normal, hence $\sigma(p^2)\equiv \sigma(p) \sigma(p)$.

Assume $\sigma(g_1g_2)\not\equiv \sigma(g_1)\sigma(g_2)$, i.e.\
$\sigma(g_1)\sigma(g_2)$ is not normal.
Then there exists a suffix of $\sigma(g_1)\sigma(g_2)$ that is not initially normal.
Since $\sigma(g_1)$ and $\sigma(g_2)$ are normal, we can write
$$
\sigma(g_1)\equiv w_1'v_1^{\epsilon_1}w_1,\qquad
\sigma(g_2)\equiv w_2v_2^{\epsilon_2}w_2'
$$
such that $v_1^{\epsilon_1}w_1w_2$ is initially normal
but $v_1^{\epsilon_1}w_1w_2v_2^{\epsilon_2}$ is not.
Therefore there exists $v_0^{\epsilon_0}$ such that
$$
v_0>v_1,\qquad
v_0\not\in S(v_1^{\epsilon_1}w_1w_2),\qquad
v_0\in S(v_1^{\epsilon_1}w_1w_2v_2^{\epsilon_2}).
$$
By the above claim, $v_0\DC v_1^{\epsilon_1}w_1w_2$
and $v_0\in S(v_2^{\epsilon_2})=\{v_2\}$.
Therefore,
$$
v_0\DC v_1,\qquad
v_0\DC w_1,\qquad
v_0\DC w_2,\qquad
v_0=v_2.
$$

Assume $w_2$ is the empty word.
Since $g_2$ is SD-conical,
we have $v_0=v_2=\apex(g_2)$ and
$v_0$ does not commute with any generator smaller than $v_0$.
This is a contradiction because $v_0\DC v_1$ and $v_0>v_1$.
Therefore $w_2$ is not the empty word, hence $S(w_2)$ is not empty.

Let $v_2'\in S(w_2)$.
Since $v_2=v_0\DC w_2$, both $v_2'$ and $v_0$
belong to $S(w_2v_2^{\epsilon_2})$.
Moreover, $v_2'\ne v_0$
because $v_2'\in S(w_2)$ and $v_0\DC w_2$ imply
$v_2'\in\supp(w_2)$ and $v_0\not\in\supp(w_2)$, respectively.
Since $S(w_2v_2^{\epsilon_2})\subset S(g_2)$, we have $|S(g_2)|\ge 2$,
which is a contradiction because $g_2$ is conical.

\smallskip
(ii)\ \
We first show that $g_1g_2$ is conical.
Assume $|S(g_1g_2)|\ge 2$.
Let $v_1=\apex(g_1)$. Then $v_1\in S(g_1)\subset S(g_1g_2)$.
Because $|S(g_1g_2)|\ge 2$,
there exists $v_2\in S(g_1g_2)\setminus\{v_1\}$.
Since $g_1$ is conical and $v_2\ne v_1$, we have $v_2\not\in S(g_1)$.
By the claim, we have $v_2\DC g_1$.
This contradicts that $g_1$ is strongly non-split.
Therefore $g_1g_2$ is conical,
and hence non-split.

Since $S(g_1)\subset S(g_1g_2)$ and $|S(g_1g_2)|=1$, we have $S(g_1g_2)=S(g_1)$,
hence $\apex(g_1g_2)=\apex(g_1)$.
Because $g_1$ is SD-conical, $\apex(g_1g_2)=\apex(g_1)$ does not
commute with any generator smaller than itself.
Therefore $g_1g_2$ is SD-conical.

The element $g_1g_2$ is strongly non-split
because $\supp(g_1g_2)\supset \supp(g_1)$ and $g_1$ is strongly non-split.
\end{proof}

\long\def\temp{
Assume that $\sigma(g_1g_2)\not\equiv \sigma(g_1)\sigma(g_2)$, i.e,
the word $\sigma(g_1)\sigma(g_2)$ is not normal.
There is a suffix of $\sigma(g_1)\sigma(g_2)$ which is not initially normal.
Since $\sigma(g_1)$ and $\sigma(g_2)$ are normal, we can write
$$
\sigma(g_1)=u_1 v_1^{\epsilon_1} w_1,\qquad
\sigma(g_2)=w_2 v_2^{\epsilon_2} u_2
$$
such that $v_1^{\epsilon_1}w_1w_2$ is initially normal,
but $v_1^{\epsilon_1}w_1w_2v_2v_2^{\epsilon_2}$ is not.
}

\begin{proposition}\label{prop:nf}
For $n\ge 1$, if $ag^n b$ is geodesic and $g$ is strongly non-split and SD-conical,
then $\sigma(ag^nb)\equiv \sigma(a)\sigma(g)^{n-1}\sigma(gb)
\equiv \sigma(a)\underbrace{\sigma(g)\cdots\sigma(g)}_{n-1}\sigma(gb)$.
\end{proposition}

\begin{proof}
Because $g$ is strongly non-split and SD-conical,
by Lemma~\ref{lem:nf}(ii),
$g^k b=g\cdot g^{k-1}b$ is strongly non-split and SD-conical
for all $1\le k\le n$.
By Lemma~\ref{lem:nf}(i), $\sigma(ag^nb)\equiv\sigma(a)\sigma(g^nb)$
and $\sigma(g^kb)\equiv\sigma(g\cdot g^{k-1}b)\equiv \sigma(g)\sigma(g^{k-1}b)$ for
$2\le k\le n$.
By an induction on $k$, we can show
$\sigma(g^kb)\equiv \underbrace{\sigma(g)\cdots\sigma(g)}_{k-1}\sigma(gb)$
for all $1\le k\le n$. Therefore
$$\sigma(ag^nb)\equiv \sigma(a)\sigma(g^nb)
\equiv \sigma(a)\underbrace{\sigma(g)\cdots\sigma(g)}_{n-1}\sigma(gb).$$
\vskip -\baselineskip
\end{proof}

\begin{lemma}\label{lem:p2sp}
Let $g_1, g_2\in G(\Gamma)$ be  strongly non-split.
Then there exists a linear order on $V(\Gamma)$
together with vertices $v_1\in\supp(g_1)$ and $v_2\in\supp(g_2)$ such that
both the $v_1$-conical conjugate of $g_1$ and the $v_2$-conical conjugate of $g_2$ are
strongly non-split and SD-conical.
\end{lemma}

\begin{proof}
First, suppose $\supp(g_1) \cap \supp(g_2) \neq\emptyset$.
Choose any $v_0\in \supp(g_1) \cap \supp(g_2)$.
Give a linear order on $V(\Gamma)$ such that
 $v_0$ is the smallest element of $V(\Gamma)$.
Let $p_1$ and $p_2$ be the $v_0$-conical conjugates of $g_1$ and $g_2$, respectively.
Then $\apex(p_1)=\apex(p_2)=v_0$ is the smallest element of $V(\Gamma)$,
hence $p_1$ and $p_2$ are SD-conical.

Now suppose $\supp(g_1) \cap \supp(g_2) = \emptyset$.
Since $g_1$ is strongly non-split, there exists $(v_1, v_2) \in \supp(g_1) \times \supp(g_2)$
such that\ $\{ v_1, v_2\}\in E(\Gamma)$.
Give a linear order on $V(\Gamma)$ such that\ $v_1 = \min V(\Gamma)$ and
$v_2 = \min \(V(\Gamma)-\{ v_1\}\)$.
For $i=1,2$, let $p_i$ be the $v_i$-conical conjugate of $g_i$.
Then $p_1$ is SD-conical because $\apex(p_1)=v_1$ is the smallest element of $V(\Gamma)$.
And $p_2$ is SD-conical because $\apex(p_2)=v_2$ is smaller than any generator other than $v_1$
and $\{v_1,v_2\}\in E(\Gamma)$.

In both cases above, $p_1$ and $p_2$ are strongly non-split
because $\supp(p_i)=\supp(g_i)$ for $i=1,2$.
\end{proof}

\section{Uniqueness of quasi-roots}

In this section, we prove Theorem~\ref{thm:AA} in three steps.

\subsection{Step 1: Quasi-roots of words}
Let $X$ be a set of symbols.
In this subsection, a letter means an element of $X^{\pm 1}$
and a word means a word on $X^{\pm1}$.
In Proposition~\ref{prop:root},
we will prove the uniqueness of primitive quasi-roots up to conjugacy
for words.

\begin{lemma}\label{lem:per}
Let $\{x_n\}_{n\ge 1}$ and $\{y_n\}_{n\ge 1}$ be sequences of letters
of period $p\ge 1$ and $q\ge 1$, respectively, such that
the first $p+q$ terms coincide, i.e.
$$\begin{array}{ll}
x_{n+p}=x_n & \mbox{for $n\ge 1$},\\
y_{n+q}=y_n & \mbox{for $n\ge 1$},\\
x_n=y_n     & \mbox{for $1\le n\le p+q$}.
\end{array}
$$
Then $x_n=y_n$ for all $n\ge 1$
and the sequence $\{x_n\}_{n\ge 1}$ has period $\gcd(p,q)$.
\end{lemma}

\begin{proof}
For $1\le n\le p$, we have
$$x_n=y_n=y_{n+q}=x_{n+q}.$$
Because $\{x_n\}_{n\ge 1}$ has period $p$, the identity $x_n=x_{n+q}$ holds for all $n\ge 1$.
Therefore $\{x_n\}_{n\ge 1}$ has period $q$, hence
it has period $\gcd(p,q)$.
By the same argument $\{y_n\}$ has period $\gcd(p,q)$.
Because $\{x_n\}$ and $\{y_n\}$ are sequences of period $\gcd(p,q)$
whose first $p+q$ terms are identical,
we have $x_n=y_n$ for all $n\ge 1$.
\end{proof}

\begin{corollary}\label{cor:root}
Let $w_1\equiv x_1\cdots x_p$ and $w_2\equiv y_1\cdots y_q$
be two primitive words of lengths $p$ and $q$, respectively.
Suppose that there exist powers $w_1^{m_1}$ and $w_2^{m_2}$
such that $w_1^2$ is a prefix of $w_2^{m_2}$ and $w_2^2$
is a prefix of $w_1^{m_1}$.
Then $p=q$ and $w_1\equiv w_2$.
\end{corollary}

\begin{proof}
Consider the infinite words
\begin{align*}
w_1^*&\equiv w_1w_1\cdots\equiv x_1\cdots x_p x_1\cdots x_p x_1\cdots x_p\cdots,\\
w_2^*&\equiv w_2w_2\cdots\equiv y_1\cdots y_q y_1\cdots y_q y_1\cdots y_q\cdots.
\end{align*}
From the hypothesis, the first $p+q$ letters of $w_1^*$ and $w_2^*$ coincide.

Let $d=\gcd(p,q)$.
By Lemma~\ref{lem:per}, $w_1^*$ and $w_2^*$ coincide and both have period $d$.
Namely, $u=x_1\cdots x_d$ is a common root of $w_1$ and $w_2$, i.e.\
$w_1\equiv u^a$ and $w_2\equiv u^b$, where $a=p/d$  and $b=q/d$.
Since $w_1$ and $w_2$ are primitive, we have $a=b=1$, hence $p=q$ and
$w_1\equiv w_2$.
\end{proof}

\begin{proposition}\label{prop:root}
Let $w$ be a word with the following two decompositions
$$
w\equiv a_1w_1^{m_1}b_1\equiv a_2w_2^{m_2}b_2,
$$
where $m_i\ge 2$ and $w_i,a_i, b_i$ are words for $i=1,2$.
Suppose that $w_1$ and $w_2$ are primitive
and that there exist non-negative constants $A$ and $ B $
such that
$$
|a_i|\le A ,\quad |b_i|\le  B ,\quad
|w|-(A  + B )\ge 2|w_i|
$$
for $i=1,2$.
Then $w_1$ and $w_2$ are cyclically conjugate.
Furthermore,  $a_1w_1a_1^{-1}=a_2w_2a_2^{-1}$ and $b_1^{-1}w_1b_1=b_2^{-1}w_2b_2$
when the words are regarded as elements of a free group.
\end{proposition}

\begin{proof}
Without loss of generality, we may assume  $|a_1|\ge |a_2|$.
Then we have a decomposition $w_2\equiv w_2'w_2''$ such that
$ a_1\equiv a_2w_2^k w_2'$ for some $k\ge 0$.
(Notice that $|a_1|+|b_2|+2|w_2|\le A  +  B  +2|w_2| \le |w| = |a_2|+|b_2|+m_2 |w_2|$
and hence $|a_1| \le |a_2| +(m_2-2) |w_2|$.)
See Figure~\ref{fig:a}.
Therefore
$$
w\equiv a_2w_2^{m_2}b_2
\equiv a_2w_2^k w_2' (w_2''w_2')^{m_2-k-1}w_2''b_2
\equiv a_1 (w_2''w_2')^{m_2-k-1}w_2''b_2.
$$

\begin{figure}
$$
\begin{xy}/r1mm/:
(0,0)="H",
"H"+(0,0)="b"; "b"+(120,0) **@{-},
"H"+(0,6)="b"; "b"+(120,0) **@{-},
"H"+(  0,0)="a"; "a"+(0,6) **@{-}, "a"+(18,0) *+!D{a_1},
"H"+( 35,0)="a"; "a"+(0,6) **@{-}, "a"+(5,0) *+!D{w_1},
"H"+( 45,0)="a"; "a"+(0,6) **@{-}, "a"+(5,0) *+!D{w_1},
"H"+( 55,0)="a"; "a"+(0,6) **@{-}, "a"+(10,0) *+!D{\cdots},
"H"+( 75,0)="a"; "a"+(0,6) **@{-}, "a"+(5,0) *+!D{w_1},
"H"+( 85,0)="a"; "a"+(0,6) **@{-}, "a"+(5,0) *+!D{\cdots},
"H"+( 95,0)="a"; "a"+(0,6) **@{-}, "a"+(12,0) *+!D{b_1},
"H"+(120,0)="a"; "a"+(0,6) **@{-},
(0,-13)="H",
"H"+(0,0)="b"; "b"+(120,0) **@{-},
"H"+(0,6)="b"; "b"+(120,0) **@{-},
"H"+(  0,0)="a"; "a"+(0,6) **@{-}, "a"+(10,0) *+!D{a_2},
"H"+( 15,0)="a"; "a"+(0,6) **@{-}, "a"+(7.5,0) *+!D{w_2^k},
"H"+( 30,0)="a"; "a"+(0,6) **@{-}, "a"+(2.5,0) *+!D{w_2'},  "a"+(7.5,0) *+!D{w_2''},
"H"+( 40,0)="a"; "a"+(0,6) **@{-}, "a"+(5,0) *+!D{w_2},
"H"+( 50,0)="a"; "a"+(0,6) **@{-}, "a"+(10,0) *+!D{\cdots},
"H"+( 70,0)="a"; "a"+(0,6) **@{-}, "a"+(5,0) *+!D{w_2},
"H"+( 80,0)="a"; "a"+(0,6) **@{-}, "a"+(5,0) *+!D{w_2},
"H"+( 90,0)="a"; "a"+(0,6) **@{-}, "a"+(5,0) *+!D{\cdots},
"H"+(100,0)="a"; "a"+(0,6) **@{-}, "a"+(10,0) *+!D{b_2},
"H"+(120,0)="a"; "a"+(0,6) **@{-},
(35,13)="a", "a"; "a"+(0,-38) **@{.},
(83,13)="a", "a"; "a"+(0,-38) **@{.},
(35,-21) *++!R{\mbox{word of length $|a_1|$}},
(83,-21) *++!L{\mbox{word of length $\lfloor B \rfloor $}},
(59,-21) *++{\mbox{word $z$}},
\end{xy}
$$
\caption{The word $w\equiv a_1w_1^{m_1}b_1\equiv a_2w_2^{m_2}b_2$}
\label{fig:a}
\end{figure}
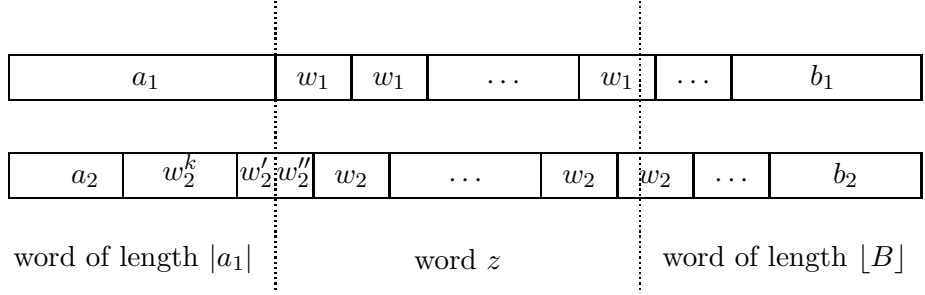

Let $z$ be the subword of $w$ obtained by removing the first $|a_1|$ letters
and the last $\lfloor B\rfloor $ letters from $w$.
From the hypothesis that $|w|-(A  +  B  )\ge 2|w_i|$, we have
$|z|\ge 2|w_i|$ for $i=1, 2$.
Notice that $z$ is a prefix of both $w_1^{m_1}$ and $(w_2''w_2')^{m_2-k-1}w_2''$.

Since $z$ is a prefix of $w_1^{m_1}$, it is of the form
$$
z\equiv w_1\cdots w_1w_1'\equiv w_1^{n_1}w_1'
$$
for some $n_1\ge 2$ and some  prefix $w_1'$ of $w_1$.
On the other hand,
since $z$ is prefix of $(w_2''w_2')^{m_2-k-1}w_2''$, it is of the form
$$
z \equiv (w_2''w_2')^{n_2}w_2'''
$$
for some $n_2\ge 2$ and some prefix $w_2'''$ of $w_2''w_2'$.

Notice that $w_2''w_2'$ is primitive because it is a cyclic conjugate of
the primitive word $w_2$.
The words $w_1$ and $w_2''w_2'$ now satisfy the hypothesis of Corollary~\ref{cor:root}.
Therefore $w_1\equiv w_2''w_2'$, hence
$w_1$ and $w_2$ are cyclically conjugate.

When the words are regarded as elements of a free group,
\begin{align*}
a_1w_1a_1^{-1}
&=(a_2w_2^kw_2')(w_2''w_2') (a_2w_2^kw_2')^{-1}
=a_2w_2^k w_2'w_2'' w_2^{-k}a_2^{-1}\\
&=a_2w_2^k w_2 w_2^{-k}a_2^{-1}
=a_2w_2a_2^{-1}.
\end{align*}
Using the same argument, we obtain
$b_1^{-1}w_1b_1=b_2^{-1}w_2b_2$.
\end{proof}

\bigskip

\subsection{Step 2: SD-conical quasi-roots}

\begin{proposition}\label{prop:2}
Let $\Gamma$ be a connected graph, and let
$$h=a_1p_1^{n_1}b_1=a_2p_2^{n_2}b_2$$
be two geodesic decompositions of\/ $h\in G(\Gamma)$,
where $n_i\ge 3$ and $a_i,b_i,p_i\in G(\Gamma)$ for $i=1, 2$.
Suppose that $p_1$ and $p_2$ are strongly non-split, SD-conical and primitive
and that there exist non-negative constants $A $ and $ B $ such that
$$
|a_i|\le A ,\quad |b_i|\le  B ,\quad
|h|-(A  +  B )\ge 3|p_i|
$$
for $i=1, 2$. Then the CGW-normal forms $\sigma(p_1)$ and $\sigma(p_2)$
are cyclically conjugate.
Furthermore, $a_1p_1a_1^{-1}=a_2p_2a_2^{-1}$ and $b_1^{-1}p_1b_1=b_2^{-1}p_2b_2$.
\end{proposition}

\begin{proof}
Without loss of generality, we may assume $|p_1|\ge |p_2|$.

Since $a_1 p_1^{n_1} b_1$ and $a_2 p_2^{n_2} b_2$ are geodesic decompositions
and since $p_1$ and $p_2$ are strongly non-split and SD-conical,
the following holds by Proposition~\ref{prop:nf}.
\begin{equation}\label{E:pro}
\sigma(h)
\equiv \sigma(a_1)\sigma(p_1)^{n_1-1}\sigma(p_1b_1)
\equiv \sigma(a_2)\sigma(p_2)^{n_2-1}\sigma(p_2b_2)
\end{equation}
Let $A '=A $ and $ B '= B +|p_1|$. Then
$$|h|-(A ' +  B ')=|h|-(A  +  B )-|p_1|\ge 3|p_1|-|p_1|=2|p_1|\ge 2|p_2|.$$
Applying Proposition~\ref{prop:root} to \eqref{E:pro}, we get the desired result.
\end{proof}

\subsection{Step 3: Proof of Theorem~\ref{thm:AA}}

\begin{proposition}\label{prop:3}
Let $\Gamma$ be a connected graph, and let
\begin{equation}\label{E:p1}
h=a_1g_1^{n_1}b_1=a_2g_2^{n_2}b_2
\end{equation}
be two geodesic decompositions of $h\in G(\Gamma)$,
where $n_i\ge 1$ and $a_i,b_i,g_i\in G(\Gamma)$ for $i=1, 2$.
Suppose that $g_1$ and $g_2$ are strongly non-split and primitive
and that there exist non-negative constants
$A$ and $B$ such that
$$
|a_i|\le A,\quad
|b_i|\le B,\quad
|h| - (A+B)\ge (2|V(\Gamma)|+1)|g_i|
$$
for $i=1,2$.
Then $g_1$ and $g_2$ are conjugate such that
$a_1g_1a_1^{-1}=a_2g_2a_2^{-1}$ and $b_1^{-1}g_1b_1=b_2^{-1}g_2b_2$.
\end{proposition}

\begin{proof}
Let $V=|V(\Gamma)|$.
Since $n_i|g_i|=|g_i^{n_i}|=|h|-(|a_i|+|b_i|)\ge |h|-(A+B)\ge (2V+1)|g_i|$,
we have $n_i\ge 2V+1\ge 2$.
Since $a_ig_i^{n_i}b_i$ is geodesic,
$g_i^{n_i}$ is also geodesic, hence
$g_i$ is cyclically reduced by Lemma~\ref{lem:c-red1}.

By Lemma~\ref{lem:p2sp},
we can choose a linear order on $V(\Gamma)$
together with vertices $v_1\in\supp(g_1)$ and $v_2\in\supp(g_2)$
such that the $v_i$-conical conjugate $p_i$ of $g_i$
is strongly non-split and SD-conical for $i=1,2$.
Note that $|g_i|=|p_i|$.

By Proposition~\ref{prop:pc},
we have the following geodesic decompositions
\begin{equation}\label{E:p2}
g_1^{n_1}=c_1 p_1^{n_1-k_1} d_1,\quad g_2^{n_2}=c_2 p_2^{n_2-k_2} d_2
\end{equation}
for some $0\le k_1, k_2 \le V-1$,
where the elements $c_i, d_i\in G(\Gamma)$ are such that
$g_i^{k_i}=c_i d_i$ is geodesic and
$g_i = c_i p_i c_i^{-1} = d_i^{-1} p_i d_i$
for $i=1,2$.
Since $g_i^{k_i}=c_i d_i$ is geodesic, we have
$$
|c_i|,~|d_i|\le | g_i^{k_i} | \le {k_i}|g_i|\le (V-1)|g_i|
$$
for $i=1,2$.
Combining \eqref{E:p1} and \eqref{E:p2},
we have the following two decompositions of $h$:
\begin{equation}\label{E:p3}
h=a_1c_1 p_1^{n_1-k_1} d_1 b_1  =a_2 c_2 p_2^{n_2-k_2} d_2 b_2.
\end{equation}
They are geodesic decomposions because both $h=a_i g_i^{n_i}b_i$ and
$g_i^{n_i}=c_i p_i^{n_i-k_i} d_i$ are geodesic.

Let $r=\max\{|g_1|, |g_2|\}$, $A'=A+(V-1)r$ and $B'=B+(V-1)r$.
Then for $i=1,2$
\begin{align*}
|a_ic_i| &\le |a_i|+|c_i|\le A+(V-1)|g_i|\le A',\\
|d_ib_i| &\le |d_i|+|b_i|\le B+(V-1)|g_i|\le B',
\end{align*}
\begin{align*}
|h|-(A'+B')
&=|h|-(A+B)-2(V-1)r\\
&\ge (2V+1)r-2(V-1)r= 3r\\
&\ge 3|g_i|=3|p_i|.
\end{align*}
Since $p_i$ is conjugate to the primitive element $g_i$, it is also primitive.
Since $n_i\ge 2V+1$ and $k_i\le V-1$, one has $n_i-k_i\ge V+2\ge 3$.
Applying Proposition~\ref{prop:2} to \eqref{E:p3}, we get
\begin{equation*}
a_1 c_1 p_1 c_1^{-1} a_1^{-1} = a_2 c_2 p_2 c_2^{-1} a_2^{-1},\qquad
b_1^{-1} d_1^{-1} p_1 d_1 b_1 = b_2^{-1} d_2^{-1} p_2 d_2 b_2.
\end{equation*}
From the identity $g_i = c_i p_i c_i^{-1} = d_i^{-1} p_i d_i$,
we obtain
$$
a_1g_1a_1^{-1}=a_2g_2a_2^{-1},\qquad
b_1^{-1}g_1b_1=b_2^{-1}g_2b_2.
$$

\end{proof}

\begin{proof}[Proof of Theorem~\ref{thm:AA}]
We are given $0\le \lambda < \frac12$, $N\ge \frac{2|V(\Gamma)|+1}{1-2\lambda}$ and two geodesic decompositions of $h$
\begin{equation}\label{E:main1}
h=a_1g_1^{n_1}b_1=a_2g_2^{n_2}b_2,
\end{equation}
where $n_i\ge N$, $a_i,b_i,g_i\in G(\Gamma)$,
$|a_i|,|b_i|\le \lambda|h|$, and
$g_i$ is strongly non-split and primitive for $i=1,2$.

Since $|h|\ge |g_i^{n_i}| = n_i |g_i| \ge N |g_i|$,
one has $|g_i|\le |h|/N$ for $i=1,2$.
Let $A=B=\lambda |h|$.
Then $|a_i|\le A$, $|b_i|\le B$ and
\begin{align*}
&|h|-(A+B)=(1-2\lambda)|h|
\ge (2 |V(\Gamma)|+1)\frac{|h|}N\ge (2 |V(\Gamma)|+1)|g_i|
\end{align*}
for $i=1,2$. Applying Proposition~\ref{prop:3},
we get the desired results.
\end{proof}

\section*{Acknowledgements}
The first author was partially supported by  NRF-2018R1D1A1B07043291.
The second author was partially supported by NRF-2018R1D1A1B07043268.

\end{document}